\documentclass{amsart}
\usepackage{amssymb}
\newcommand{\cC}{{\mathcal C}}

\newcommand{\be}{{\bf 1}}

\newcommand{\BZ}{{\mathbb Z}}

\newcommand{\BQ}{{\mathbb Q}}
\newcommand{\BC}{{\mathbb C}}

\newcommand{\Rep}{\mbox{Rep}}
\newcommand{\FP}{{\rm FPdim}}

\newcommand{\Id}{\mbox{Id}}
\newcommand{\End}{\mbox{End}}

\title{Pre-modular categories of rank 3}
\author{Victor Ostrik}
\email{vostrik@darkwing.uoregon.edu}
\address{Department of Mathematics, 1222 University of Oregon, Eugene OR 97403-1222}
\thanks{The author was partially supported by NSF grant DMS-0098830 and 
DMS-0111298.}
\date{March 2005}
\begin{document}
\begin{abstract} We classify ribbon semisimple monoidal categories with 
three isomorphism classes of simple objects over the field of complex numbers.
\end{abstract} 
\maketitle
\section{Introduction}
Let $k$ be an algebraically closed field of characteristic 0.
A fusion category $\cC$ over $k$ is a $k-$linear semisimple rigid monoidal
category with finitely many simple objects and finite dimensional spaces
of morphisms, such that the endomorphism algebra of the neutral object is $k$,
see \cite{ENO}. The rank of a fusion category is just the number of 
isomorphism classes of simple objects. In \cite{O2} we classified all fusion
categories of rank 2. A similar classification problem for rank 3 seems
to be out of reach at the moment. For example it is not known whether the
number of fusion categories of rank 3 is finite. In this note we classify
the fusion categories $\cC$ of rank 3 under an additional assumption that
$\cC$ admits a ribbon structure, see \cite{BK}. Recall that the ribbon 
fusion categories are called pre-modular, see \cite{Mu}.

{\bf Main Theorem.} There are exactly 7 fusion categories of rank 3 admitting
a structure of ribbon category.

The proof of this theorem is based on the properties of $S-$matrix 
(see \cite{BK}) and is an exercise in Galois theory. This note was
inspired by \cite{GK} where the authors classified 
the fusion rings of modular tensor categories with small number of
simple objects and small fusion coefficients via computer search. 
It would be very interesting to answer the following

{\bf Question.} Is it true that there are only finitely many ribbon 
categories of a given finite rank?

This question is a special case of question in \cite{O2}. On the other hand the positive
answer to this question would imply the conjecture by Z.~Wang that there are just
finitely many modular tensor categories of a given rank. 

It is interesting to observe that contrary to the case of rank 2 not all fusion categories of rank 3 admit a ribbon structure, see Remark in section 4.5.

After this note was finished D.~Nikshych showed to the author reference \cite{CP}
where the fusion rings of modular tensor categories with three simple objects were classified
under some unitarity assumptions. Also I was informed
by Z.~Wang that all modular tensor categories of rank $\le 4$ are now classified, see
\cite{W}. Still we hope that our treatment will be useful for some readers.

This note was written while the author enjoyed the hospitality of the Institute for
Advanced Study. I am happy to thank this institution. I am grateful to Dmitri Nikshych and
Zhenghan Wang for very useful comments.

\section{Recollections}
In this section we review necessary facts on pre-modular categories.

\subsection{Dimension} Recall (see e.g. \cite{ENO}) that a fusion category 
$\cC$ is pivotal if it is endowed with a functorial tensor isomorphisms 
$M\to M^{**}$ for any $M\in \cC$. In a pivotal fusion category $\cC$ one
defines for any object $M\in \cC$ its dimension $d_M\in k$, see \cite{BK}.
We have the following properties, see {\em loc. cit.}:

(a) $d_M$ defines a homomorphism $d_M: K(\cC)\to k$;

(b) Assume that $M$ is simple object. Then $d_M\ne 0$.

\subsection{$\tilde S-$matrix} Let $\cC$ be a ribbon category (see \cite{BK}
for a definition). Recall that in a ribbon category the balancing isomorphism
$\theta \in \End (\Id_{\cC})$ is defined. For a simple object $X\in \cC$ let
$\theta_X$ denote the scalar by which $\theta$ acts on $X$. Vafa's theorem
(see \cite{V,BK,Ev}) states that 

(a) the numbers $\theta_X$ are roots of unity.  

Let $\{ V_i\}_{i\in I}$ be a set of (reprsentatives of isomorphism classes of)
simple objects in $\cC$ and let $i\mapsto i^*$ be a unique involution of $I$
such that $(V_i)^*\cong V_{i^*}$. Let $V_i\otimes V_j=
\sum_{k\in I}N_{ij}^kV_k$ and let $\theta_i:=\theta_{V_i}, d_i:=d_{V_i}$. 
Define the matrix $\tilde S$ by the formula (see \cite{BK}):
$$\tilde S_{ij}=\theta_i^{-1}\theta_j^{-1}\sum_{k\in I}N_{i^*j}^k\theta_kd_k.$$

We have the following properties, see \cite{BK}:

(b) the matrix $\tilde S$ is symmetric $\tilde S_{ij}=\tilde S_{ji}$;

(c) For any $i\in I$ the assignment $\phi_i(V_j)=\tilde S_{ij}/d_i$ defines
a homomorphism of rings $K(\cC)\to k$.

Recall (see \cite{BK}) that a ribbon category $\cC$ is called modular if
the matrix $\tilde S$ is non-degenerate.

\section{Ribbon based rings of rank 3}
\subsection{} Let $k,l,m,n$ be nonnegative integers subject to the condition
$$k^2+l^2=lm+kn+1. \eqno(*)$$ 
Let $K(k,l,m,n)$ be the based ring with the basis
$1,X,Y$ and the multiplication given by
$$X^2=1+mX+kY,\; Y^2=1+lX+nY,\; XY=YX=kX+lY.$$
The following Proposition gives the classification of the based rings of
rank 3, see \cite{ENO}, Example in section 8.10.

{\bf Proposition.} Let $K$ be a based ring of rank 3. Then  
either $K=K(\Rep(\BZ/3\BZ))$ or $K=K(k,l,m,n)$.

Note that we have an obvious isomorphism of the based rings 
$K(k,l,m,n)=K(l,k,n,m)$. 

\subsection{Symmetric categories} Recall that a ribbon category $\cC$ is 
called symmetric if the square of the braiding is the identity. Equivalently,
the $\tilde S-$matrix of the category $\cC$ has rank 1 (see e.g. \cite{Mu}).
It is proved by Deligne \cite{De} that for any symmetric fusion category $\cC$
there exists a finite group $G$ and an equivalence $\cC \simeq \Rep(G)$. 

{\bf Lemma.} Let $G$ be a finite group with 3 irreducible representations.
Then either $G=\BZ/3\BZ$ or $G=S_3$.

{\bf Proof.} The Landau estimate (see \cite{L,O2}) gives $|G|\le 6$. The rest
is easy.  

{\bf Corollary.} Let $\cC$ be a symmetric category of rank 3. Then either
$K(\cC)=K(\Rep(\BZ/3\BZ))$ or $K(\cC)=K(0,1,0,1)$.

\subsection{Non-modular and non-symmetric categories} Assume that the
category $\cC$ is not symmetric and is not modular. It follows from \cite{Mu}
Corollary 2.16 that the category $\cC$ has a non-trivial symmetric 
subcategory. This
subcategory has 2 simple objects (say $\be$ and $X$) and thus is equivalent
to $\Rep(\BZ/2\BZ)$. Hence $K(\cC)=K(0,1,0,n)$. Observe that $d_X=1$ since
otherwise $d_Y=0$. Thus $\tilde S-$matrix looks like
$$\tilde S=\left(\begin{array}{ccc}1&1&d_Y\\1&\theta_X^{-2}&\theta_X^{-1}d_Y\\
d_Y&\theta_X^{-1}d_Y&\theta_Y^{-2}(1+\theta_X+n\theta_Yd_Y)
\end{array}\right) .$$
 
Since the second column should give a homomorphism $K(\cC)\to \BC$ we get
$\theta_X=1$ (except, possibly, the case $n=0$). Since the third column is a 
$d_Y$ times a homomorphism $K(\cC)\to \BC$
distinct from $d$, we get $\theta_Y^{-2}(1+\theta_X+n\theta_Yd_Y)=
d_Y\bar d_Y$ where $\bar d_Y$ is a root of the equation $y^2=2+ny$ distinct
from $d_Y$. Thus we have $\theta_Y^{-2}(1+\theta_X+n\theta_Yd_Y)=-2$ or
equivalently $nd_Y=-2(\theta_Y+\theta_Y^{-1})$. Assume that $n>1$. Then
$d_Y$ is irrational and after applying a Galois avtomorphism to the last
equation we have $ny_+=-2(\theta +\theta^{-1})$ where $y_+$ is the positive
root of the equation $y^2=2+ny$ and $\theta$ is a some root of unity. But
note that $y_+^2>ny_+$ and hence $y_+>n$ and $ny_+>n^2$. On the other hand
clearly $|\theta +\theta^{-1}|\le 2$. Thus $n^2<4$ and we get the 
contradiction. Thus we have proved

{\bf Proposition.} Assume that $\cC$ is nor symmetric neither modular. Then
$K(\cC)=K(0,1,0,n)$ where $n=0,1$.

{\bf Remark.} It is reasonable to expect that if $K(\cC)=K(0,1,0,n)$ for
some fusion category $\cC$ then $n\le 2$ (see Remark in Section 4.5). 
But unfortunately we don't know
how to prove that $n$ is bounded by any constant.

\subsection{Modular categories} In this section we assume that $\cC$ is a
modular category such that $K(\cC )=K(k,l,m,n)$. Let $\phi_1,\phi_2,\phi_3$ be
the three distinct homomorphisms $K(\cC)\to \BC$; we assume that $\phi_1$
coincides with the dimension function and denote $\phi_i(X)=x_i, 
\phi_i(Y)=y_i$ for $i=1,2,3$. We can assume that $\tilde S-$matrix looks like
$$\tilde S=\left( \begin{array}{ccc}1&x_1&y_1\\x_1&x_1x_2&y_1x_3\\y_1&x_1y_2&
y_1y_3\end{array}\right).$$
Since $\tilde S-$matrix is symmetric we have $x_1y_2=y_1x_3$. It is easy to
see that $x_1y_2=y_1x_3\ne 0$.

The absolute Galois group $Gal(\bar \BQ/\BQ)$ acts on the set 
$\{ \phi_1,\phi_2,\phi_3\}$. Thus we have a homomorphism $Gal(\bar \BQ/\BQ)\to
S_3$. Let us denote the image of this homomorphism by $G$. It is known 
\cite{dBG} that the group $G$ is abelian (see also \cite{CG}, \cite{ENO} 
Appendix). Thus we have 3 possibilities:
$G$ is trivial, $G=\BZ/3\BZ$ and $G=\BZ/2\BZ$.

{\bf Case 1.} $G$ is trivial. Then all numbers $x_i, y_i$ are rational and
hence integer. In particular $\FP(X)$ and $\FP(Y)$ are integers. Then again
the Landau estimate (see \cite{ENO} 8.38) gives $\FP(\cC)\le 6$ and the only
possibility is $K(\cC)=K(0,1,0,1)$.

{\bf Case 2.} $G=\BZ/3\BZ$. The group $G$ permutes homomorphisms $\phi_1,
\phi_2,\phi_3$ cyclically. Thus applying the elements of $G$ to the identity
$x_1y_2=y_1x_3$ we get new identities $x_2y_3=y_2x_1$ and $x_3y_1=y_3x_2$.
Equivalently $x_1y_2=x_2y_3=x_3y_1=:\lambda$. Recall that
$\lambda \ne 0$. Thus $(x_1,x_2,x_3)=\lambda (y_2^{-1},y_3^{-1},y_1^{-1})$.
The numbers $x_i$ are the roots of the polynomial
$x^3-(m+l)x^2+(ml-k^2-1)x+l$ (the characteristic polynomial of the operator
of multiplication by $X$ in $K(k,l,m,n)$) and the numbers $y_i$ are the
roots of the polynomial $y^3-(n+k)y^2+(nk-l^2-1)y+k$. The Vieta Theorem
implies:
$$-l=x_1x_2x_3=\frac{\lambda^3}{y_1y_2y_3}=\lambda^3/(-k) \Rightarrow 
\lambda^3=lk; \eqno (1)$$
$$m+l=x_1+x_2+x_3=\lambda (y_1^{-1}+y_2^{-1}+y_3^{-1})=\lambda 
\frac{nk-l^2-1}{(-k)}; \eqno (2)$$
$$ml-k^2-1=x_1x_2+x_2x_3+x_3x_1=\lambda^2\frac{y_1+y_2+y_3}{y_1y_2y_3}=
\lambda^2\frac{n+k}{(-k)}. \eqno (3)$$
Now equation (2) implies that $\lambda$ is rational except, possibly,
the case $m+l=0$ (in the latter case $K(k,l,m,n)=K(1,0,0,0)$). Equation
(1) then says that $\lambda >0$ and equations (2), (3) imply
$nk-l^2-1<0$ and $ml-k^2-1<0$. But we know from (*) that 
$nk-l^2-1+ml-k^2-1=-3$ and hence we can assume that $nk-l^2-1=-1$ and
$ml-k^2-1=-2$. It is easy to see that these equations imply $k=1$ and
$l=1,2$. But the case $l=2$ is impossible since then $\lambda =\sqrt[3]{kl}$
is irrational. Thus the only possibility is $K(\cC)=K(1,1,1,0)$.
Thus we found that in case 2 we have 2 possibilities: either 
$K(\cC)=K(1,1,1,0)$ or $K(k,l,m,n)=K(1,0,0,0)$.

{\bf Case 3.} $G=\BZ/2\BZ$. In this case there are two subcases:

(a) $G$ fixes $\phi_1$ and permutes $\phi_2$ and $\phi_3$. The identity 
$x_1y_2=y_1x_3$ implies $x_1y_3=y_1x_2$ and hence $x_2y_2=x_3y_3$. Thus
$\phi_2(X\otimes Y)=\phi_3(X\otimes Y)$. Since $\phi_2\ne \phi_3$ we see
that $(XY)^2$ should lie in the subspace of $K(\cC)$ spanned by $1$ and $XY$. 
Now $(XY)^2=(kX+lY)^2=k^2X^2+2klXY+l^2Y^2=k^2(mX+kY)+l^2(lX+nY)\mod{<1,XY>} =
(k^2m+l^3)X+(k^3+l^2n)Y\mod{<1,XY>}$. This vector should be proportional to
$XY=kX+lY$, hence $(k^2m+l^3)l=(k^3+l^2n)k$. We see that if $p$ is a prime 
divisor of $k$ then $p$ divides $l$ and the relation (*) then shows that
$p$ divides $1$. Thus $k\le 1$ and similarly $l\le 1$. Thus in this case
we have that either $K(\cC)=K(0,1,0,n)$ or $K(\cC)=K(1,1,1,0)$. In the first
case we have $x_1=1$ (otherwise $y_1=0$) and $y_1$ is a root of equation
$y^2=2+ny$. This equation has a rational root only for $n=1$. Thus we have
2 possibilities $K(\cC)=K(0,1,0,1)$ and $K(\cC)=K(1,1,1,0)$.

(b) $G$ does not fix $\phi_1$. This is most difficult case. We can assume 
that $G$ permutes $\phi_1$ and $\phi_2$ and fixes $\phi_3$. Thus the identity 
$x_1y_2=y_1x_3$ implies $x_2y_1=y_2x_3$ and hence $x_1x_2=x_3^2$. Thus by 
Vieta Theorem $x_3^3=x_1x_2x_3=-l$. Set $t:=x_3$ and $s:=y_3$; then $s$ and 
$t$ are integers. Assuming $s\ne 0$ we have
$$k=s+t^2s,\; l=-t^3,\; m=t-ts^2-\frac{s^2+1}{t},\; n=s+\frac{t^4-1}{s}.$$
Also
$$x_1x_2=\frac{-l}{x_3}=t^2,\; x_1+x_2=m+l-x_3=-ts^2-t^3-\frac{s^2+1}{t},$$
$$y_1y_2=\frac{-k}{y_3}=-(t^2+1),\; y_1+y_2=n+k-y_3=t^2s+s+\frac{t^4-1}{s}.$$
The equation $x_1y_2=y_1x_3$ gives $x_1=y_1^2x_3/(y_1y_2)=
-\frac{t}{t^2+1}y_1^2$ and, similarly, $x_2=-\frac{t}{t^2+1}y_2^2$.
Thus we have $x_1+x_2=-\frac{t}{t^2+1}(y_1^2+y_2^2)$ or, equivalently,
$$-ts^2-t^3-\frac{s^2+1}{t}=-\frac{t}{t^2+1}((s(t^2+1)+\frac{t^4-1}{s})^2+
2(t^2+1)).$$
After simple transformations we have
$$\frac{s^2+1}{t^2}+t^2-s^2t^2-2t^4=\frac{(t^2-1)^2}{s^2}(t^2+1)$$
and thus
$$\frac{s^2}{t^2}(1-t^4)+\frac{1}{t^2}+t^2-2t^4=\frac{(1-t^2)^2}{s^2}(t^2+1).$$
After dividing by $1-t^2$ we get
$$\frac{s^2}{t^2}(1+t^2)+\frac{1+t^2+2t^4}{t^2}=\frac{1-t^2}{s^2}(1+t^2)$$
or, equivalently,
$$\frac{s^2}{t^2}+\frac{1}{t^2}+\frac{2t^2}{t^2+1}+\frac{t^2}{s^2}=
\frac{1}{s^2}.$$
But this is impossible since the LHS is greater than $\frac{s^2}{t^2}+
\frac{t^2}{s^2}\ge 2$ and the RHS is $\le 1$. 

Thus we have two possibilities: either $s=0$ or $t^2=1$. Assume first that
$t^2=1$. Then $t=-1$, and $k=2s, l=1, m=2s^2, n=s$. It is not difficult to
check that in this case the $\tilde S$ matrix is symmetric. We have
$-y_1=y_1x_3=\tilde S_{23}=\frac{1}{\theta_X\theta_Y}
(2s\theta_Xx_1+\theta_Yy_1)$. Recall that $x_1=-\frac{t}{t^2+1}y_1^2=
\frac{1}{2}y_1^2$. Hence $sy_1=-\theta_Y-\frac{\theta_Y}{\theta_X}$. After
applying the Galois automorphism we can assume that $y_1>0$ and we have
an inequality $sy_1\le 2$. On the other hand $y_1$ is a root of the polynomial
$y^2-2sy-2$ and hence $y_1^2=2sy_1+2>2sy_1 \Rightarrow y_1>2s$. Thus we get
$2s^2<1$ and hence $s=0$.

Consider now the case $s=0$. Then we have $k=0, l=1, m=0$. The 
$\tilde S-$matrix looks like 
$$\tilde S=\left( \begin{array}{ccc}1&1&y_1\\1&1&-y_1\\y_1&y_2&0\end{array}
\right) $$
where $y_1,y_2$ are the roots of the equation $y^2=2+ny$. Since $\tilde S$ is
symmetric we have $y_2=-y_1$ and hence $n=0$. Thus $K(\cC)=K(0,1,0,0)$.

Summarizing we can state

{\bf Proposition.} Assume that a fusion category $\cC$ of rank 3 admits a
structure of modular category. Then we have the following possibilities
for $K(\cC)$: $K(\Rep(\BZ/3\BZ))$, $K(1,0,0,0)$, $K(0,1,0,1)$, $K(1,1,1,0)$.

\subsection{List of possible based rings} We have proved

{\bf Theorem.} Assume that a fusion category of rank 3 admits a ribbon
structure. Then we have the following possibilities for $K(\cC)$:
$K(\Rep(\BZ/3\BZ))$, $K(0,1,0,0)$, $K(0,1,0,1)$, $K(1,1,1,0)$.

\section{Identification of tensor categories}
In this section we describe all fusion categories with Grothendieck rings 
given by Theorem 3.5.

\subsection{} $K(\cC)=K(\Rep(\BZ/3\BZ)$. In this case possible fusion 
categories are classified by $H^3(\BZ/3\BZ,k^*)=\BZ/3\BZ$ (note that 
$\mbox{Aut}(\BZ/3\BZ)$ acts trivially on this cohomology group, see e.g.
\cite{EGO}). Thus there are 3 such categories. But only the category with
trivial associativity constraint admits a structure of ribbon category,
see e.g. \cite{Qu}. This structure is not unique: we can have symmetric
category $\cC =\Rep(\BZ/3\BZ)$ and modular category $\Rep (\widehat{sl}(3)_1)$.

\subsection{} $K(\cC)=K(0,1,0,0)$. The fusion ring $K(0,1,0,0)$ is well
known in conformal field theory, it represents the fusion rules of the
Ising model. We have an isomorphism $K(0,1,0,0)=K(\Rep(\widehat{sl}(2)_2)$.
Thus acoording to \cite{Ke} (see also \cite{KW}) there are two fusion
categories $\cC$ such that $K(\cC)=K(0,1,0,0)$. One of them is 
$\Rep(\widehat{sl}(2)_2)$ and the second can be obtained from the first one
by applying some Galois automorphism; also both categories can be constructed
using the quantum group $U_q(sl(2))$ for $q=\sqrt[8]{1}$, see \cite{BK}.

\subsection{} $K(\cC)=K(1,1,1,0)$. Observe that $K(1,1,1,0)\boxtimes 
K(\Rep(\BZ/2\BZ))=K(\Rep(\widehat{sl}(2)_5)$. Thus it follows from \cite{Ke} 
(see also \cite{KW}) that there are exactly three fusion categories $\cC$ 
with such Grothendieck ring; one category is a subcategory 
$\Rep(\widehat{so}(3)_5)$  of representations with integer spin in 
$\Rep(\widehat{sl}(2)_5)$; two others are Galois conjugate to this one.
Thus all three categories admit a ribbon structure. Also all three categories
can be realized using the quantum group $U_q(sl(2))$ for $q=\sqrt[7]{1}$,
see \cite{BK}.

\subsection{} $K(\cC)=K(0,1,0,1)$. Observe that $K(0,1,0,1)=K(\Rep(S_3))$.
It was established by T.~Chmutova that there are 3 fusion categories with
such Grothendieck ring, see \cite{EGO}. It is easy to see that 
$\tilde S-$matrix for such category necessarily has rank 1 and hence any
braided structure on $\cC$ is symmetric. Thus only $\cC =\Rep(S_3)$ of these 
3 categories has a ribbon structure.

\subsection{} Summarizing the results of the previous sections we can state
the main result of this note.

{\bf Main Theorem.} There are exactly 7 fusion categories of rank 3 admitting
a ribbon structure: $\Rep(\BZ/3\BZ)$, $\Rep(S_3)$, $\Rep(\widehat{sl}(2)_2)$,
$\Rep(\widehat{so}(3)_5)$ and the Galois conjugates of two latter categories.

{\bf Remark.} It is interesting to note that there exists a fusion category of rank 3
(and hence with commutative Grothendieck ring) which does not admit a ribbon
structure. Namely let $\cC$ be the fusion category attached to the affine $sl_2$ on
level 10 and let $A\in \cC$ be the commutative $\cC-$algebra of type $E_6$, see \cite{KO}. 
Then the category $\Rep A$ of right $A-$modules has a structure of fusion category,
see {\em loc. cit.} This fusion category contains a tensor subcategory corresponding
to the ends of long legs and the triple vertex of the graph $E_6$ which is fusion
category of rank 3 with Grothendieck ring $K(0,1,0,2)$. It follows from the Main Theorem
above that this category does not admit a ribbon structure.

It seems reasonable to expect that the Grothendieck ring of any
fusion category of rank 3 is either isomorphic to $K(0,1,0,2)$ or listed in Theorem 3.5.

\end{document}